\documentclass[12pt]{article}
\usepackage{amsmath,amssymb,latexsym}
\usepackage{theorem}
\usepackage{graphicx}
\newtheorem{thm}{Theorem}
\newtheorem{cor}[thm]{Corollary}

{\theorembodyfont{\rm}   }
{\theorembodyfont{\rm}   }
{\theorembodyfont{\rm}   }
{\theorembodyfont{\rm}   }

\def \lcm{\mbox{lcm}} \def \gcd{\mbox{gcd}}

\renewcommand{\labelenumi}{\setlength{\labelwidth}{\leftmargin}
   \addtolength{\labelwidth}{-\labelsep}
   \hbox to \labelwidth{\theenumi.\hfill}}

\begin{document}
\title{The last chapter of the Disquisitiones of Gauss}

\author{Laura Anderson, Jasbir S. Chahal, Jaap Top}
\date{ }




\maketitle

\begin{abstract}
This exposition reviews what exactly Gauss 
asserted and what did he prove in the last chapter of {\sl Disquisitiones Arithmeticae} about
dividing the circle into a given number of equal parts. In other words, what did Gauss claim and actually prove concerning the roots of unity and the construction
of a regular polygon with a given number of sides.
Some history of Gauss's solution is briefly recalled, and in particular 
many relevant classical references are provided which we believe deserve to be better known.
\end{abstract}

\section{Introduction}
A remark by J-P.~Serre made towards one of us concerning what the textbook \cite{jasbir} asserts about Gauss' results on the construction
of regular $n$-gons, provided the initial motivation to write this note.
A second source of inspiration came from
D.~Surya Ramana's expository paper \cite{Ramana} on Carl Friedrich Gauss (1777--1855). Various useful additions as well as interesting historical
notes are presented in a slightly later text by B.~Sury \cite{Sury}.
In 1801 Gauss published his {\it Disquisitiones Arithmeticae} \cite{gauss}.   The purpose of the present article is to elaborate on the remark of Serre and the comments by Ramana and Sury concerning the last (seventh) chapter of this celebrated textbook.
Its last chapter is devoted to the study of roots of unity, i.e., the (complex) solutions to an equation $x^n=1$. 
 Gauss applied the results he obtained on roots of unity to the problem of constructing  a regular polygon of $n$ sides
 or, what amounts to the same, divide the circle into $n$
 equal parts. Here (and everywhere in the present paper) `constructing' means constructing by using a straightedge and compass  alone. 

Ramana writes (p.~64 of \cite{Ramana})
\begin{quote}
$\ldots$ {\sl Gauss proved that a regular polygon with $n$ sides, where $n$ is an odd integer, can be
constructed using a ruler and compass alone if and only if $n$ can be expressed as a product of distinct Fermat primes.}
\end{quote}

Gauss indeed discovered and formulated a general result, namely a necessary and sufficient condition on positive integer $n$ for a regular polygon with $n$ sides to be constructible.  He states it in the very last Article~366 of Chapter~VII of the
{\it Disquisitiones} as follows (whenever we quote from  Gauss's text, the 1986 revised English translation \cite{gauss} is used) :
\begin{center}
\scalebox{.6}{\includegraphics{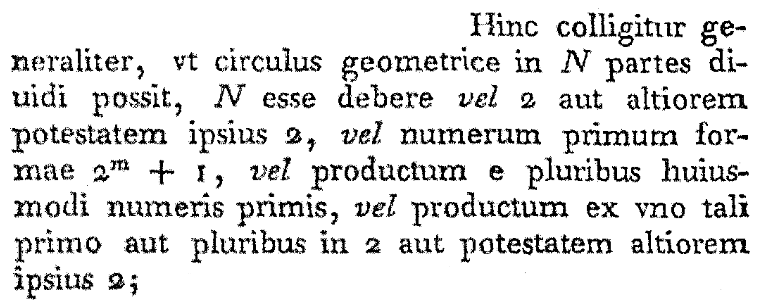}}
\label{Gauss theorem}
\end{center}
{``In general therefore in order to be able to divide the circle geometrically into $N$ parts, $N$ must be 2 or a higher power of 2, \emph{or} a prime number of the form $2^m+1$, \emph{or} the product of several prime numbers of this form, \emph{or} the product of one or several such primes into 2 or a higher power of 2;''}

\vspace{\baselineskip}
However, Gauss also
writes (Article~365) that ``{\sl the limits  of the present work exclude $\ldots$ demonstration here.}'' and (Article~335)
``{\sl $\ldots$ although we could discuss them {\rm (he refers to roots of unity in ${\mathbb C}$ here)}
in all their generality, we reduce them to the simplest case $\ldots$, both for the sake of brevity and in order that the new
principles of this theory may be more easily understood.}'' With  this last sentence he apparently
restricts himself to $p$-th roots of unity, for (odd) prime numbers $p$. 

The precise modern formulation of the statement by Gauss quoted above is
 \begin{thm}\label{thm1} A regular $N$-gon is constructible if and only if $N>1$ is of the form
\[
N=2^r p_1 p_2 \ldots p_s
\]
with $r,s\in{\mathbb Z}_{\geq 0}$ and $p_j$  distinct Fermat primes, i.e., primes of the form $2^m+1$.
 \end{thm}

Article~365 of the {\sl Disquisitiones} contains a proof of the fact that if $2^m+1$ is prime, then
$m=2^\nu$, as well as the observation that Fermat studied these numbers and that Euler discovered $2^{32}+1$ to be divisible by $641$.
 Although Fermat believed {\em all} integers $2^{2^\nu}+1$ to be prime, it is an old conjecture (for example,
see \cite[p.~315]{EPascal} that they are prime for only finitely many
integers $\nu$. Likely, $\nu=4$ is the largest of these.

Gauss finishes the {\sl Disquisitiones} by listing all integers below $300$ satisfying the
criterion of Theorem~\ref{thm1}:
\begin{center}
\scalebox{.4}{\includegraphics{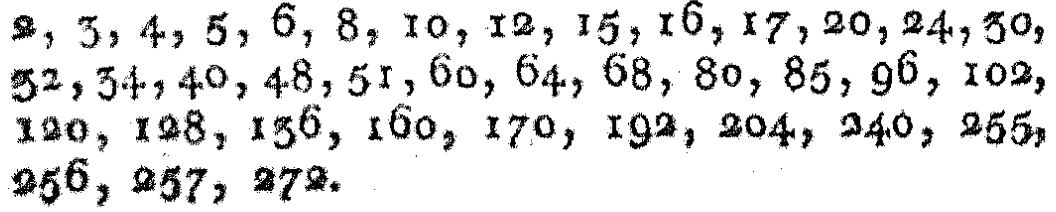}}
\end{center}
Many sources assert that Gauss proved Theorem~\ref{thm1}, for example Chapter~5 of the history textbook \cite{jasbir}
written by the second author of the current paper. After receiving a copy  of these history notes,
Jean-Pierre Serre pointed out to him (as already alluded to in the first lines of this paper) that the {\sl Disquisitiones} does not contain a complete proof of the result.
In the remainder of the present note we review some of the famous history of the problem of constructing regular $n$-gons,
and its current mathematical formulation. Note, incidentally,
that Serre's remark was observed by others long before him: see, for
example, the title of Pierpont's 1895 text \cite{Pierpont}, a paper
discussing among other things a proof of Theorem~\ref{thm1}.
Also Sury's text \cite[p.~48]{Sury} points out that this theorem is
stated but not proven in the {\sl Disquisitiones}.
\subsection{Constructions with straightedge and compass}

Constructions rely on two axioms
already presented in Euclid's {\sl Elements} (Book~1):
\begin{enumerate}
\item[1)]Any two distinct points may be joined by one and only one (straight) line.
\item[2)] In the plane, a circle may be drawn with any point as its center and any finite length segment as its radius.
\end{enumerate}

\noindent Clearly, they are equivalent to allowing the use of two basic tools:
\begin{enumerate}
\item[1)]
A straightedge to join two distinct points, and
\item[2)] a compass to draw circles with already determined centers and radii.
 \end{enumerate}

A convenient approach to this identifies the plane with the complex plane, i.e., with the field of complex numbers
${\mathbb C}$. This is explained in the first chapter of Klein's 1895 textbook \cite{klein}, and in terms of ${\mathbb R}^2$ in \S37 of Hilbert's textbook \cite{Hilbert} of 1899, and 
using field theory and ${\mathbb C}$ in \S65 of Van der Waerden's 1930 classic \cite{vdW}.
It leads to the general criteria below.
\begin{thm}\label{constructible}
A complex number $z$ is constructible starting from $0$ and $1$
 if and only if one of the following equivalent statements holds:
 \begin{itemize}
     \item[\text{\rm (i)}] 
 Field extensions $
{\mathbb Q}=K_0\subseteq K_1\subseteq \ldots \subseteq K_e\subset{\mathbb C}
$
exist such that every $K_j\subseteq K_{j+1}$ is quadratic, and $z\in K_e$.
\item[\text{\rm (ii)}] The splitting field of ${\mathbb Q}(z)$ over $\mathbb{Q}$ has degree a power of $2$.
\end{itemize}
\end{thm}

The Greeks were  interested in knowing what kind of figures are possible to construct.  In Book~1 (Proposition~1) of his {\sl Elements}, Euclid shows how to construct an equilateral triangle on a given line segment $AB$. 
In Book~1 (Proposition~46) it is shown how to construct a square.  Finally, Book~4 (Proposition~11) describes a method  to construct a regular pentagon.

A regular $n$-gon, or simply  $n$-gon, is an $n$ sided convex figure in the Euclidean plane with all sides equal, or equivalently with all the interior angles equal to each other.  For example, an equilateral triangle is a 3-gon and a regular pentagon is a 5-gon.  
Since it is possible to bisect any given angle, one can obtain a $2n$-gon from an $n$-gon.  

More generally, given two positive integers $n,m$, constructing the regular $n$-gon
and the regular $m$-gon is equivalent to constructing the regular
$\lcm(n,m)$-gon. Namely, put $k:=\lcm(n,m)$. Constructing a regular $k$-gon
is the same as constructing an angle of size $2\pi/k$. Having this,
evidently all integral multiples of this angle can be constructed,
so in particular those of size $2\pi/n$ and $2\pi/m$.
Vice versa, from these two angles one can also construct angles of
size $2\pi x/n+2\pi y/m$ for arbitrary integers $x,y$.
Taking $x,y$ such that $xm+yn=\gcd(n,m)$ (such integers $x,y$ exist
by the extended Euclidean algorithm), one obtains
\[
x\frac{2\pi}{n}+y\frac{2\pi}{m}=\frac{2\pi}{nm}(xm+yn)=\frac{2\pi}{\lcm(n,m)}.
\]
This is classical, and briefly recalled in Article~336 of the {\sl Disquisitiones}.

The construction of an $n$-gon is equivalent to the construction of the numbers $e^{2\pi im/n}=\cos(2\pi m/n)+i\sin(2\pi m/n)$ in the
complex plane, given the numbers $0$ and $1$. In other words, one tries to construct
from $0$ and $1$ the complex solutions of $x^n=1$.

Ramana's exposition \cite{Ramana} amply discusses the breakthrough by Gauss
concerning the construction of $n$-gons.  In the words of Gauss himself (see~\cite{arch2}),
 \begin{quote}
``During the winter of 1796 (my first semester in
G\"ottingen), I had already discovered everything related to the
separation of the roots of the equation
 $$\frac {x^p-1}{x-1}=0$$
into \emph{two} groups, on which the beautiful theorem\ldots
depends\ldots.  After intensive consideration of the relation of all
the roots to one another on arithmetical grounds, I succeeded during
a holiday in Braunschweig, on the morning of the day alluded to
[March 29, 1796] (before I had got out of bed), in viewing this
relation in the clearest way, so that I could immediately make
special application to the 17-side and to the numerical
verification\ldots. I announced this discovery\ldots in May or June 1796.''
 \end{quote}

In Article~365 of the \emph{Disquisitiones} Gauss  observes,
 \begin{quote}
``It is certainly astonishing that although the geometric divisibility
of the circle into three and five parts was already known in
Euclid's time, nothing was added to this discovery for 2000 years.
And all geometers had asserted that, except for those sections and
the ones that derive directly from them (that is, division into 15,
$3 \cdot 2^\mu$, $5 \cdot 2^\mu$, and $2^\mu$ parts), there
are no others that can be effected by geometric constructions.''
 \end{quote}

According to the German mathematician Heinrich Weber (1842--1913) famous for the
Kronecker-Weber theorem, Gauss requested that the $17$-gon
 be engraved on his tombstone (footnote in \S~97 of the first edition of \cite{Weber} :
\begin{center}
\scalebox{.38}{\includegraphics{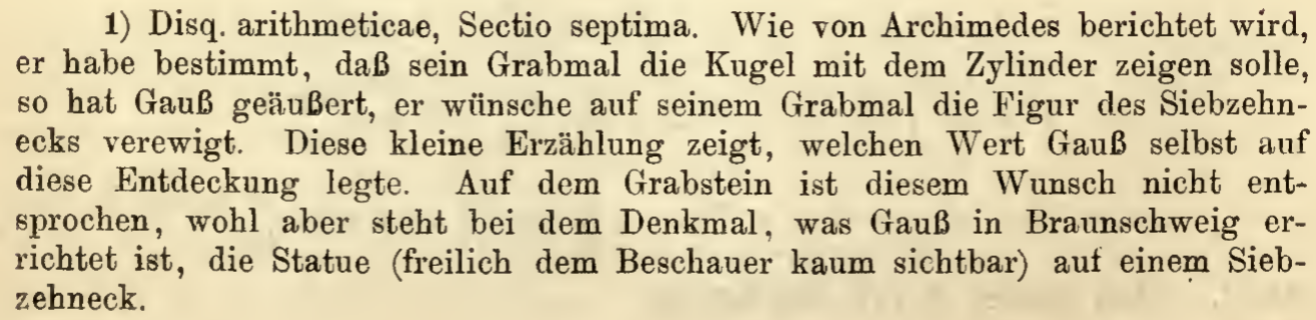}} 
\end{center}
Translation: 1) Disq. arithmeticae, Chapter VII. As is recorded that Archime\-des ascertained that his epitaph
show a ball with a cylinder, Gauss expressed his wish that his epitaph immortalizes the figure of the $17$-gon.
This little tale shows the importance Gauss himself put on his discovery. On the epitaph this wish was not
fulfilled, yet on the monument erected for Gauss in his native town Braunschweig, the statue is placed on a (barely visible
for the visitor) $17$-gon.

\vspace{\baselineskip}
Both Ramana \cite{Ramana} and Sury \cite{Sury} present a detailed derivation of Gauss' formula for $\cos(2\pi/17)$; whereas
Ramana's paper does not contain the resulting formula itself,
Sury's one does.
Note that the formula given in Article~365 of the {\sl Disquisitiones}
actually contains a typo: it reads
\begin{center}
\scalebox{.4}{\includegraphics{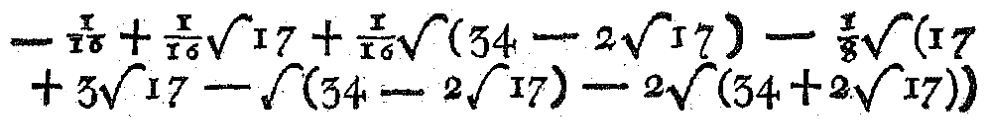}} 
\end{center}
and the rightmost minus sign on the first line should be a $+$, as
confirmed by (e.g.) Sury's short and clear derivation \cite[bottom of p.~46]{Sury} (using his notation, $\alpha_3/2=\cos(2\pi/17)$)
and the resulting correct formula on \cite[p.~47]{Sury}.
From the modern perspective as given in Theorem~\ref{constructible} above, Gauss found
the explicit tower of quadratic field extensions
\[
{\mathbb Q}\subseteq{\mathbb Q}(\sqrt{17})
\subseteq {\mathbb Q}\left(\sqrt{34-2\sqrt{17}}\right)
\subseteq{\mathbb Q}(\cos(2\pi/17))\subseteq 
{\mathbb Q}(e^{2\pi i/17}).
\]
Although both the formula and the series of quadratic
extensions lead to a  construction method
of the $17$-gon, an actual construction is not discussed in the {\sl Disquisitiones}.
According to R.C.~Archibald's historical survey \cite{arch2}, J.F.~Pfaff in 1802
wrote a letter to Gauss in which he quotes the construction from a letter he received from
C.F.~Pfleiderer. Archibald reproduces the part
of Pfaff's letter containing the construction:
\begin{center}
\scalebox{.48}{\includegraphics{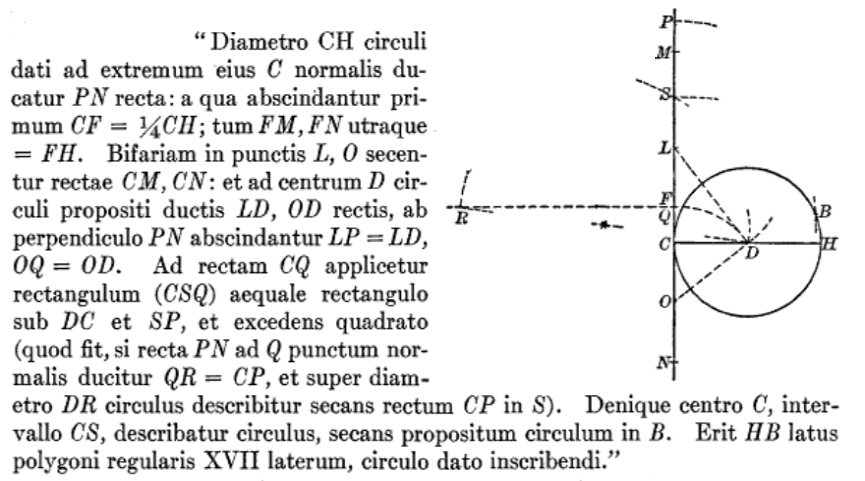}} 
{(taken from R.C.~Archibald \cite{arch2}.)}\label{Archibald}
\end{center}

The next Fermat prime after $2^4+1=17$ is $2^8+1=257$.
In 1831, Friedrich Julius Richelot in his PhD thesis \cite{Richelot} (supervised
by Carl Jacobi) presents a construction of the $257$-gon. 
This was also done in 1887 by E. Pascal \cite{Pascal257} and in 1893 in a short book 
by Hugo Schwendenwein \cite{Schw}. For the final known Fermat prime $2^{16}+1=65537$, Johann Gustav Hermes in 1894 published a construction of the
65537-gon \cite{Hermes}. 
The practical difficulty in constructing a precise 65537-gon is obvious, for even 257 is so large that on a page of a
book, the 257-gon will look like a circle.

\section{Proving Gauss's  claim}

A contemporary proof of Theorem~\ref{thm1} goes like this; compare \cite{Pierpont}, \cite[\S~19.3]{Garling} and
many similar texts.

 For  any integer $n \ge 1$, put
 \[
\omega = e^{2\pi i/n}\in {\mathbb C}.
 \]
The complex numbers
 $$1, \omega, \omega^2, \ldots, \omega^{n-1}$$
are the zeros of the polynomial
\[x^n - 1\]
of degree $n$.  

As remarked earlier, constructing a regular $n$-gon is equivalent to 
constructing regular $p^m$-gons, for every prime power $p^m$ (with $m>0$)
exactly dividing $n$.
Writing $\zeta=e^{2\pi ip^{-m}}$, by Theorem~\ref{constructible} the latter is equivalent to
the existence of field extensions
\[
{\mathbb Q}=K_0\subset K_1\subset \ldots \subset K_e={\mathbb Q}(\zeta)
\]
in which every $K_j\subset K_{j+1}$ is quadratic.

The extension ${\mathbb Q}(\zeta)/{\mathbb Q}$ has degree $\phi(p^m)=p^{m-1}(p-1)$ where $\phi$ is Euler's $\phi$-function. By assumption the
extension $K_e/K_0$ has degree $2^e$. Hence if the $p^m$-gon is
constructible then either $p=2$ (and  $m>0$ is arbitrary),
or $m=1$ and $p$ is a Fermat prime.
Vice versa, if $\phi(p^m)=2^e$, then the observation that
 ${\mathbb Q}(\zeta)/{\mathbb Q}$ is a Galois extension of degree $2^e$ with abelian Galois group, implies by the fundamental theorem of Galois theory
that $\zeta$ satisfies the condition for constructibility given in Theorem~\ref{constructible}.
This finishes the proof of Theorem~\ref{thm1}.

 \bigskip
A special instance of Theorem~\ref{thm1} is that an $18$-gon, and therefore an angle of $20^{\circ}$,
cannot be constructed whereas a $6$-gon, so an angle of $60^{\circ}$, can. This shows:
\begin{cor}\label{coroll}
Given only a compass and straightedge, it is impossible to trisect all given angles.
\end{cor}
In 1837, Pierre Laurent Wantzel \cite{want} published a  (different) proof of Corollary~\ref{coroll}, although it is
an immediate consequence of Theorem~\ref{thm1}. An interesting discourse
why Gauss might have emphasized constructibility while somewhat neglecting the impossibility of constructions was recently
given by L\"{u}tzen \cite{L1}, \cite{L2}. These papers moreover discuss earlier proofs of Wantzel's theorem
by R.~Descartes (1637), J.\'{E}.~Montucla  (1754), and N.J.A.N.~de Condorcet (1775), all preceding the publication of the {\sl Disquisitiones}.

In what follows, we quote what Gauss himself says of the proof of Theorem~\ref{thm1} he left out.

Concerning $p^{\lambda}$-th roots of unity with $\lambda>1$ and $p$ an odd prime, Article~336 of the {\sl Disquisitiones} asserts
\begin{quote}
``to determine a polygon of $p^{\lambda}$ sides we necessarily require the solution of $\lambda-1$
equations of degree $p$. Even though the following theory could be extended to this case also,
nevertheless we could not avoid so many equations of degree $p$, and there is no way of reducing their degree if $p$
is prime.''
\end{quote}
In today's terminology, this seems to say that if $\zeta=e^{2\pi i/p}$ and $\omega=e^{2\pi i p^{-\lambda}}$,
then the extension $K(\omega)/K(\zeta)$ has degree $p^{\lambda-1}$ and in fact a chain of intermediate fields
\[
{\mathbb Q}(\zeta)=L_0\subset K_1\subset \ldots \subset K_{\lambda-1}={\mathbb Q}(\omega)
\]
exists in which each $K_{j+1}$ has degree $p$ over $K_j$.

In his discussion of polynomials such that  $\zeta=e^{2\pi i/p}$  can be expressed in terms of their zeros (with $p$ an odd prime number not of the form $1+2^m$), he writes (\emph{Disquisitiones}, Article~365), the capitalization is by Gauss himself:
 \begin{quote}
``{\sc we can show with all  rigor that these higher-degree equations cannot be avoided in any way nor can they be reduced to 
lower-degree equations.} The limits of the present work exclude this demonstration here, but we issue this warning lest anyone attempt to achieve geometric constructions for sections other than the ones suggested by our theory (e.g. sections into 7, 11, 13, 19, etc. parts) and so spend his time uselessly.''
 \end{quote}

In modern terms, this is interpreted as the claim that the degree of ${\mathbb Q}(\zeta)$ over ${\mathbb Q}$, which equals $p-1$,
is divisible by an odd number $>1$. 
 \bigskip

\noindent For $n$  a Fermat prime, Gauss states ({\sl Disquisitiones}, Article~365; note that he uses $P$ for $2\pi$),
 \begin{quote}
``Whenever$\ldots$ the value of $n$ is 3, 5, 17, 257, 65537, etc. the
sectioning of the circle is reduced to quadratic equations only, and
the trigonometric functions of the angles $P/n$, $2P/n$, etc.
can be expressed by square roots which are more or less complicated
(according to the size of $n$).  Thus in these cases the division of
the circle into $n$ parts or the inscription of a regular polygon of
$n$ sides can be accomplished by geometric constructions.''
 \end{quote}

In todays language, this claims that the field ${\mathbb Q}(\zeta)$ with $\zeta=e^{2\pi i/p}$ and $p$ a Fermat prime,
admits a chain of subfields
\[
{\mathbb Q}=K_0\subset K_1\subset \ldots \subset K_e={\mathbb Q}(\zeta),
\]
with each $K_{j+1}$ quadratic over $K_j$.

Write $\Phi_n$ for the minimal monic polynomial with rational coefficients
having $\omega=e^{2\pi i/n}$ as a zero. The degree of $\Phi_n$ equals the degree of
${\mathbb Q}(e^{2\pi i/n})$ over ${\mathbb Q}$. Gauss determined in Article~341 of the
{\sl Disquisitiones} that for $n=p$ a prime number, $\Phi_p=x^{p-1}+x^{p-2}+\ldots +x+1$.
This may be regarded as a proof of the assertion in the second quote by Gauss mentioned in this section.
Note that Entry \#136 in Gauss's Diary \cite{Gaussdiary} written in 1808 essentially says that he was
able to prove that for every $n$ the degree of $\Phi_n$ equals $\phi(n)$.
Section~IV of Wantzel's 1837 paper \cite{want} contains the assertion that a slight adaptation
of Gauss's argument for determining $\Phi_p$ gives $\Phi_{p^r}$ for all $r>0$. Although he refrained from
explicitly observing it,
Klein's determination (1895) of $\Phi_{p^2}$ sketched below in fact extends to the case $\Phi_{p^r}=(x^{p^r}-1)/(x^{p^{r-1}}-1)$, confirming Wantzel's assertion.
Earlier in 1850 Serret \cite{Serret} in a less
elementary way computed $\Phi_n$ when $n$ is a power of a prime number.
Kronecker \cite{Kronecker} in 1854 determined $\Phi_n$ for all $n$, confirming what Gauss asserted in
his diary 46 years earlier.
A historical survey containing many interesting classical references may be found
in \cite[Kapitel~IV \S~12]{EPascal}.
In Section~III.8 of \cite{klein} Klein determines $\Phi_n$ for $n=p^2$ and $p$ prime,
namely, $\Phi_{p^2}=(x^{p^2}-1)/(x^p-1)$. \text In fact, replacing $x$ by $x+1$ allows one to apply
the well known Eisenstein criterion for irreducibility. Klein continues by  observing that  the degree of this polynomial is
a multiple of $p$. Hence for $p$ odd and $n$ any multiple of $p^2$, the $n$-gon cannot be constructed.
\section{Conclusion}
Gauss in the {\sl Disquisitiones} stated a correct 
and complete criterion regarding constructibility
of a regular $n$-gon (Theorem~\ref{Gauss theorem}
above). 

Moreover, by determining (stated in modern terms)
for any odd prime $p$ the minimal polynomial $\Phi_p$ of $e^{2\pi i/p}$, he shows that indeed only
the prime $2$ and Fermat primes can occur as prime factors
of $n$, in case the $n$-gon is constructible.

Gauss correctly states the reason why no odd prime
$p$ can have the property $p^2|n$ in case the
$n$-gon is constructible. However, he did not provide a proof in this case. Based on a diary
entry he wrote $7$ years after the publication of
the {\sl Disquisitiones}, and on a very simple
argument published later by Klein and probably already known to Wantzel, it seems reasonable to
assume that indeed Gauss did have a complete
proof of this, although maybe not yet in 1801.

Finally by classical arguments, to show that indeed for integers $n$
of the kind described by Gauss the $n$-gon is
constructible, it suffices to show that for
any Fermat prime $F_\nu=2^{2^\nu}+1$ the
$F_\nu$-gon is constructible. Although we have
not discussed it here, arguably Gauss knew
how to achieve this (see Articles~343-364 and
the summary he presents
in the first lines of Article~365). The
example $\nu=2$, so $F_\nu=17$ presented by him
illustrates this method.

\section*{Acknowledgement}
It is our pleasure to thank Sudhir Ghorpade for helpful advise
concerning publication of the present note. We are grateful for
the comments of a referee, in particular informing us of the
paper \cite{Sury} and of the typo in the original publication of
the {\sl Disquisitiones} which we had not noticed.
 
\bigskip

\noindent Laura Anderson, 1065 West 10210 South, South Jordan, UT 84095, USA, laurabcannon@gmail.com
 \medskip

\noindent Jasbir S. Chahal, Department of Mathematics, Brigham Young University, Provo, UT  84602, USA, jasbir@math.byu.edu
 \medskip

\noindent Jaap Top, Bernoulli Institute, University of Groningen, P. O. Box 407, 9700~AK Groningen, The Netherlands, j.top@rug.nl

\end{document}